\input amssym.def 
\input amssym.tex
\documentstyle[a4,12pt]{article} 

\input amssym.def 
\def \C{{\Bbb C}}
\def \N{{\Bbb N}}

\newtheorem{lemme}{LEMMA}[section]

\newtheorem{defi}[lemme]{DEFINITION}
\newtheorem{prop}[lemme]{PROPOSITION}
\newtheorem{cor}[lemme]{COROLLARY}

\title{Tensor products of $C(X)$-algebras over $C(X)$}
\author{Etienne Blanchard}
\date{}

\begin{document}
\maketitle 

\begin{abstract}

Given a Hausdorff compact space $X$, we study the 
\mbox{${\rm C}^*$}-(semi)-norms on the algebraic tensor 
product $A\otimes_{alg,C(X)} B$ of two $C(X)$-algebras $A$ and $B$ over $C(X)$. 
In particular, if one of the two $C(X)$-algebras defines a continuous field of 
\mbox{${\rm C}^*$}-algebras over 
$X$, there exist minimal and maximal \mbox{${\rm C}^*$}-norms on 
$A\otimes_{alg,C(X)} B$ but there does not exist any \mbox{${\rm C}^*$}-norm 
on $A\otimes_{alg,C(X)} B$ in general. 

\begin{center}
\end{center}
\end{abstract}

\setcounter{section}{-1}
\section{Introduction}
\indent\indent 
Tensor products of \mbox{${\rm C}^*$}-algebras have been extensively studied 
over the last decades (see references in \cite{la}). One of the main results 
was obtained by M. Takesaki in \cite{ta} where he proved that the spatial 
tensor product $A\otimes_{\min} B$ of two \mbox{${\rm C}^*$}-algebras $A$ and 
$B$ always defines the minimal \mbox{${\rm C}^*$}-norm on the algebraic tensor 
product $A\otimes_{alg} B$ of $A$ and $B$ over the complex field $\C$.

More recently, G.G. Kasparov constructed in \cite{kus} a tensor product over 
$C(X)$ for $C(X)$-algebras. The author was also led to introduce in \cite{bla1} 
several notions of tensor products over $C(X)$ for $C(X)$-algebras and to study 
the links between those objects. 

Notice that E. Kirchberg and 
S. Wassermann have proved in \cite{kirch} that the subcategory of continuous 
fields over a Hausdorff compact space is not closed under such tensor products 
over $C(X)$ and therefore, in order to study tensor products over $C(X)$ of 
continuous fields, it is natural to work in the $C(X)$-algebras framework.

\vskip12pt 
Let us introduce the following definition: 
\begin{defi} Given two $C(X)$-algebras $A$ and $B$, we denote by 
${\cal I} (A,B)$ the involutive ideal of the algebraic tensor product 
$A\otimes_{alg} B$ generated by the elements ${(fa)\otimes b}-a\otimes (fb)$, 
where $f\in C(X)$, $a\in A$ and $b\in B$.
\end{defi}

Our aim in the present article is to study the \mbox{${\rm C}^*$}-norms on the 
algebraic tensor product $(A\otimes_{alg} B)/{\cal I} (A,B)$ of two 
$C(X)$-algebras $A$ and $B$ over $C(X)$ and to see how one can 
enlarge the results of Takesaki to this framework.

\vskip12pt
We first define an ideal ${\cal J} (A,B)\subset A\otimes_{alg} B$ which 
contains ${\cal I} (A,B)$ such that every \mbox{${\rm C}^*$}-semi-norm on 
$A\otimes_{alg} B$ which is zero on ${\cal I} (A,B)$ is also zero on 
${\cal J} (A,B)$ and we prove that there always exist a minimal 
\mbox{${\rm C}^*$}-norm $\|\;\|_m$ and a maximal \mbox{${\rm C}^*$}-norm 
$\|\;\|_M$ on the quotient $(A\otimes_{alg} B)/{\cal J} (A,B)$.

We then study the following question of G.A. Elliott (\cite{ellio}): when do 
the two ideals ${\cal I} (A,B)$ and ${\cal J} (A,B)$ coincide?

\vskip12pt
The author would like to express his gratitude to C. Anantharaman-Delaroche and 
G. Skandalis for helpful comments. He is also very indebted to S. Wassermann 
for sending him a preliminary version of \cite{kirch} and to J. Cuntz who 
invited him to the Mathematical Institute of Heidelberg.

\section{Preliminaries}\label{john}
We briefly recall here the basic properties of $C(X)$-algebras. 

\vskip12pt
Let $X$ be a Hausdorff compact space and $C(X)$ be the 
\mbox{${\rm C}^*$}-algebra of 
continuous functions on $X$. For $x\in X$, define the morphism 
$e_x:C(X)\rightarrow\C$ of evaluation at $x$ and denote by $C_x(X)$ 
the kernel of this map.
\begin{defi}\label{lakral} (\cite{kus}) 
A $C(X)$-algebra is a \mbox{${\rm C}^*$}-algebra $A$ endowed 
with a unital morphism from $C(X)$ in the center of the multiplier 
algebra $M(A)$ of $A$. 

We associate to such an algebra the unital 
$C(X)$-algebra ${\cal A}$ generated by $A$ and $u[C(X)]$ in $M[A\oplus C(X)]$ 
where $u(g)(a\oplus f)=ga\oplus gf$ for $a\in A$ and $f,g\in C(X)$.
\end{defi}

For $x\in X$, denote by $A_x$ the quotient of $A$ by the closed ideal 
$C_x(X)A$ and by $a_x$ the image of $a\in A$ in the fibre 
$A_x$. Then, as 
\begin{center}
$\| a_x\| =\inf\{\| [1-f+f(x)]a\| ,f\in C(X)\}$, 
\end{center}
the map $x\mapsto\| a_x\|$ is upper semi-continuous for all 
$a\in A$ (\cite{ri}).

Note that the map $A\rightarrow\oplus A_x$ is a monomorphism since if 
$a\in A$, there is a pure state $\phi$ on $A$ such that $\phi (a^*a)=
\| a\|^2$. As the restriction of $\phi$ to $C(X)\subset M(A)$ 
is a character, there exists $x\in X$ such that $\phi$ factors 
through $A_x$ and so $\phi (a^*a)=\| a_x\|^2$.

Let $S(A)$ be the set of states on $A$ 
endowed with the weak topology and let ${\cal S}_X (A)$ be the subset of states 
$\varphi$ whose restriction to $C(X)\subset M(A)$ is a character, i-e such that 
there exists an $x\in X$ (denoted $x=p(\varphi )$) verifying $\varphi (f)=f(x)$ 
for all $f\in C(X)$. Then the previous paragraph implies that the set of pure 
states $P(A)$ on $A$ is included in ${\cal S}_X (A)$.

\vskip12pt 
Let us introduce the following notation: if ${\cal E}$ is a Hilbert $A$-module 
where $A$ is a \mbox{${\rm C}^*$}-algebra, we will denote 
by ${\cal L}_A({\cal E} )$ or simply ${\cal L} ({\cal E} )$ the set of bounded 
$A$-linear operators on ${\cal E}$ which admit an adjoint (\cite{kas}).
\begin{defi}\label{nike} (\cite{bla1}) Let $A$ be a $C(X)$-algebra.

A $C(X)$-representation 
of $A$ in the Hilbert $C(X)$-module ${\cal E}$ is a morphism $\pi :
A\rightarrow{\cal L} ({\cal E})$ which is $C(X)$-linear, i.e. such that for 
every $x\in X$, the representation 
$\pi_x={\pi\otimes e_x}$ in the Hilbert space ${\cal E}_x=
{\cal E}\otimes_{e_x}\C$ factors through a representation of $A_x$. 
Furthermore, if $\pi_x$ is a faithful representation of 
$A_x$ for every $x\in X$, $\pi$ is said to be a field of faithful representations of $A$.

A continuous field of states on $A$ is a $C(X)$-linear map $\varphi : 
A\rightarrow C(X)$ such that 
for any $x\in X$, the map $\varphi_x=e_x\circ\varphi$ defines a state on $A_x$.
\end{defi}

If $\pi$ is a $C(X)$-representation of the $C(X)$-algebra $A$, the 
map $x\mapsto\|\pi_x(a)\|$ is lower semi-continuous since $\langle\xi ,
\pi (a)\eta\rangle\in C(X)$ for every $\xi ,\eta\in{\cal E}$. Therefore, if $A$ 
admits a field of faithful representations $\pi$, the map 
$x\mapsto\| a_x\| =\|\pi_x(a)\|$ is continuous for every $a\in A$, which 
means that $A$ is a continuous field of \mbox{${\rm C}^*$}-algebras over $X$ 
(\cite{di}).

The converse is also true (\cite{bla1} th\'eor\`eme 3.3): 
given a separable $C(X)$-algebra $A$, the following assertions are equivalent:
\begin{enumerate}
\item $A$ is a continuous field of \mbox{${\rm C}^*$}-algebras over $X$,
\item the map $p:{\cal S}_X (A)\rightarrow X$ is open,
\item $A$ admits a field of faithful representations.
\end{enumerate}

\section{\mbox{${\rm C}^*$}-norms on 
$(A\otimes_{alg} B)/{\cal J} (A,B)$}\label{blabla} 
\begin{defi}
Given two $C(X)$-algebras $A$ and $B$, we define the involutive ideal 
${\cal J} (A,B)$ of the algebraic tensor product $A\otimes_{alg} B$ of elements 
$\alpha\in A\otimes_{alg} B$ such that $\alpha_x=0$ in 
$A_x\otimes_{alg} B_x$ for every $x\in X$.
\end{defi}

By construction, the ideal ${\cal I} (A,B)$ is included in ${\cal J} (A,B)$.

\begin{prop}\label{fraise} 
Assume that $\|\;\|_\beta$ is a \mbox{${\rm C}^*$}-semi-norm on the algebraic 
tensor product $A\otimes_{alg} B$ of two $C(X)$-algebras $A$ and $B$.

If $\|\;\|_\beta$ is zero on the ideal ${\cal I} (A,B)$, then 
\begin{center}
$\|\alpha\|_\beta =0$ for all $\alpha\in{\cal J} (A,B)$.
\end{center}
\end{prop}{\mbox{{\bf Proof: }}} 
Let $D_\beta$ be the Hausdorff completion of $A\otimes_{alg} B$ for 
$\|\;\|_\beta$. By construction, $D_\beta$ is a quotient of 
$A\otimes_{\max} B$. 
Furthermore, if $C_\Delta$ is the ideal of $C(X\times X)$ of functions which 
are zero on the diagonal, the image of $C_\Delta$ in $M(D_\beta )$ is zero. 

As a consequence, the map from $A\otimes_{\max} B$ onto $D_\beta$ factors 
through the quotient ${A{\mathop{\otimes}\limits^M}_{C(X)} B}$ of 
${A\otimes_{\max} B}$ by ${C_\Delta\times (A\otimes_{\max} B)}$. 

But an easy diagram-chasing argument shows that 
$(A{\mathop{\otimes}\limits^M}_{C(X)} B)_x={A_x\otimes_{\max} B_x}$ for every 
$x\in X$ (\cite{bla1} corollaire 3.17) and therefore the image of 
${\cal J} (A,B)\subset A\otimes_{\max} B$ in 
$A{\mathop{\otimes}\limits^M}_{C(X)} B$ is zero. $\Box$
\subsection{The maximal \mbox{${\rm C}^*$}-norm}
\begin{defi} 
Given two $C(X)$-algebras $A_1$ and $A_2$, we denote by $\|\;\|_M$ the 
\mbox{${\rm C}^*$}-semi-norm 
on $A_1\otimes_{alg} A_2$ defined for $\alpha\in A_1\otimes_{alg} A_2$ by
\begin{center} 
$\|\alpha\|_M=
\sup\{\| (\sigma_1^x\otimes_{\max}\sigma_2^x)(\alpha )\| ,x\in X\}$ 
\end{center} 
where $\sigma_i^x$ is the map $A_i\rightarrow (A_i)_x$.
\end{defi}

As $\|\;\|_M$ is zero on the ideal ${\cal J} (A_1,A_2)$, if we identify 
$\|\;\|_M$ with the \mbox{${\rm C}^*$}-semi-norm induced on 
$(A_1\otimes_{alg} A_2)/{\cal J} (A_1,A_2)$, we get: 
\begin{prop}
The semi-norm $\|\;\|_M$ is the maximal \mbox{${\rm C}^*$}-norm on the 
quotient $(A_1\otimes_{alg} A_2)/{\cal J} (A_1,A_2)$.
\end{prop}{\mbox{{\bf Proof: }}} 
By construction, $\|\;\|_M$ defines a \mbox{${\rm C}^*$}-norm on 
$(A_1\otimes_{alg} A_2)/{\cal J} (A_1,A_2)$. 
Moreover, as the quotient ${A_1{\mathop{\otimes}\limits^M}_{C(X)} A_2}$ of 
${A_1\otimes_{\max} A_2}$ by 
${C_\Delta\times (A_1\otimes_{\max} A_2)}$ maps injectively in 
\begin{center}
$\mathop{\oplus}\limits_{x\in X} (A_1{\mathop{\otimes}\limits^M}_{C(X)} A_2)_x=
\mathop{\oplus}\limits_{x\in X} (\, (A_1)_x\otimes_{\max} (A_2)_x)$, 
\end{center}
the norm of $A_1{\mathop{\otimes}\limits^M}_{C(X)} A_2$ coincides on the dense 
subalgebra $(A_1\otimes_{alg} A_2)/{\cal J} (A_1,A_2)$ with $\|\;\|_M$. But we 
saw in proposition \ref{fraise} that if $\|\;\|_\beta$ is a 
\mbox{${\rm C}^*$}-norm on the algebra $(A\otimes_{alg} B)/{\cal J} (A,B)$, 
the completion of $(A\otimes_{alg} B)/{\cal J} (A,B)$ for $\|\;\|_\beta$ is a 
quotient of $A{\mathop{\otimes}\limits^M}_{C(X)} B$. $\Box$ 

\subsection{The minimal \mbox{${\rm C}^*$}-norm}
\begin{defi} Given two $C(X)$-algebras $A_1$ and $A_2$, we define the semi-norm 
$\|\;\|_m$ on $A_1\otimes_{alg} A_2$ by the formula 
\begin{center}
$\|\alpha\|_m=
\sup\{\| (\sigma_1^x\otimes_{\min}\sigma_2^x)(\alpha )\| ,x\in X\}$ 
\end{center} 
where $\sigma_i^x$ is the map $A_i\rightarrow (A_i)_x$ and we denote by 
$A_1{\mathop{\otimes}\limits^m}_{C(X)} A_2$ the Hausdorff completion of 
$A_1\otimes_{alg} A_2$ for that semi-norm.
\end{defi}
{\sl Remark:} In general, the canonical map 
$(A_1{\mathop{\otimes}\limits^m}_{C(X)} A_2)_x\rightarrow 
(A_1)_x\otimes_{\min} (A_2)_x$ is not a monomorphism (\cite{kirch}).

\vskip12pt 
By construction, $\|\;\|_m$ induces a \mbox{${\rm C}^*$}-norm on 
$(A_1\otimes_{alg} A_2)/{\cal J} (A_1,A_2)$. 
We are going to prove that this \mbox{${\rm C}^*$}-norm defines the minimal 
\mbox{${\rm C}^*$}-norm on the 
involutive algebra $(A_1\otimes_{alg} A_2)/{\cal J} (A_1,A_2)$. 

\vskip12pt
Let us introduce some notation.

Given two unital $C(X)$-algebras $A_1$ and $A_2$, let $P(A_i)\subset 
{\cal S}_X(A_i)$ denote 
the set of pure states on $A_i$ and let 
$P(A_1)\times_XP(A_2)$ denote the closed subset of $P(A_1)\times P(A_2)$ of 
couples $(\omega_1,\omega_2)$ such that $p(\omega_1)=p(\omega_2)$, where 
$p:P(A_i)\rightarrow X$ is the restriction to $P(A_i)$ of the map 
$p:{\cal S}_X(A_i)\rightarrow X$ defined in section \ref{john}. 

\begin{lemme}\label{elise} 
Assume that $\|\;\|_\beta$ is a \mbox{${\rm C}^*$}-semi-norm on the algebraic 
tensor product $A_1\otimes_{alg} A_2$ of two unital $C(X)$-algebras $A_1$ and 
$A_2$ which is zero on the ideal ${\cal J} (A_1,A_2)$ and define the 
closed subset $S_\beta\subset{P(A_1)\times_XP(A_2)}$ of couples 
$(\omega_1,\omega_2)$ such that 
\begin{center}
${\Bigl| (\omega_1\otimes\omega_2)(\alpha )\Bigr|}\leq 
\|\alpha\|_\beta$ for all $\alpha\in A_1\otimes_{alg} A_2$. 
\end{center}

If $S_\beta \not = P(A_1)\times_XP(A_2)$, there exist self-adjoint 
elements $a_i\in A_i$ such that $a_1\otimes a_2\not\in{\cal J} (A_1,A_2)$ but 
$(\omega_1\otimes\omega_2)(a_1\otimes a_2)=0$ for all 
couples $(\omega_1,\omega_2)\in S_\beta$. 
\end{lemme}{\mbox{{\bf Proof: }}} 
Define for $i=1,2$ the adjoint action $ad$ of the unitary group ${\cal U}(A_i)$ 
of $A_i$ on the pure states space $P(A_i)$ by the formula
\begin{center}
$[(ad_u)\omega ](a)=\omega (u^*au)$. 
\end{center}
Then $S_\beta$ is invariant under the product action $ad\times ad$ of 
${\cal U}(A_1)\times {\cal U}(A_2)$ and we can therefore find non empty open 
subsets $U_i\subset P(A_i)$ which are invariant under the action of ${\cal U}(A_i)$ 
such that $(U_1\times U_2)\cap S_\beta =\emptyset$.

Now, if $K_i$ is the complement of $U_i$ in $P(A_i)$, the set 
\begin{center}
$K_i^\perp =\{ a\in A_i$ $|$ $\omega (a)=0$ for all $\omega\in K_i\}$
\end{center}
is a non empty ideal of $A_i$ and furthermore, if $\omega\in P(A_i)$ is zero on 
$K_i^\perp$, then $\omega$ belongs to $K_i$ (\cite{gli} lemma 8,\cite{ta}).

As a consequence, if $(\varphi_1,\varphi_2)$ is a point of $U_1\times_XU_2$, 
there exist non zero self-adjoint elements $a_i\in K_i^\perp$ such that 
$\varphi_i(a_i)=1$. If $x=p(\varphi_i)$, this implies in particular that 
$(a_1)_x\otimes (a_2)_x\not =0$, and hence 
$a_1\otimes a_2\not\in{\cal J} (A_1,A_2)$. $\Box$

\begin{lemme} (\cite{ta} theorem 1) Let $A_1$ and $A_2$ be two unital 
$C(X)$-algebras. 

If the algebra $A_1$ is an abelian algebra, there exists only one 
\mbox{${\rm C}^*$}-norm on the quotient 
$(A_1\otimes_{alg} A_2)/{\cal J} (A_1,A_2)$.
\end{lemme}{\mbox{{\bf Proof: }}} 
Let $\|\;\|_\beta$ be a \mbox{${\rm C}^*$}-semi-norm on $A_1\otimes_{alg} A_2$ 
such that for all $\alpha\in A_1\otimes_{alg} A_2$, $\|\alpha\|_\beta =0$ if 
and only if $\alpha\in{\cal J} (A_1,A_2)$.

\vskip12pt
If $\rho\in P(A_\beta )$ is a pure state on the Hausdorff 
completion $A_\beta$ of $A_1\otimes_{alg} A_2$ for the semi-norm 
$\|\;\|_\beta$, then for every $a_1\otimes a_2\in A_1\otimes_{alg} A_2$, 
\begin{center} 
$\rho (a_1\otimes a_2)=\rho (a_1\otimes 1)\rho (1\otimes a_2)$ 
\end{center}
since $A_1\otimes 1$ is included in the center of $M(A_\beta )$. Moreover, if 
we define the states $\omega_1$ and $\omega_2$ by the formulas $\omega_1 (a_1)=
\rho (a_1\otimes 1)$ and $\omega_2(a_2)=\rho (1\otimes a_2)$, then $\omega_2$ 
is pure since $\rho$ is pure, 
and $(\omega_1,\omega_2)\in P(A_1)\times_XP(A_2)$. It follows that 
$P(A_\beta )$ is isomorphic to $S_\beta$. 

In particular, if $a_1\otimes a_2\in A_1\otimes_{alg} A_2$ verifies 
\begin{center}
$(\omega_1\otimes\omega_2)(a_1\otimes a_2)=0$ for all couples 
$(\omega_1,\omega_2)\in S_\beta$,
\end{center}
the element $a_1\otimes a_2$ is zero in $A_\beta$ and therefore belongs to the 
ideal ${\cal J} (A_1,A_2)$. Accordingly, the previous lemma implies that 
$P(A_1)\times_XP(A_2)=S_\beta =P(A_\beta )$. 

\vskip12pt 
As a consequence, we get for every $\alpha\in A_1\otimes_{alg} A_2$ 
\begin{center}
\begin{tabular}{ll}
$\|\alpha\|_\beta^2$&$=\sup\{\rho (\alpha^*\alpha ), 
\rho\in P(A_\beta )\}$\\
& $=\sup\{ (\omega_1\otimes\omega_2)(\alpha^*\alpha ), 
(\omega_1,\omega_2)\in P(A_1)\times_XP(A_2)\}$
\end{tabular}
\end{center}
But that last expression does not depend on $\|\;\|_\beta$, and hence the 
unicity. $\Box$

\begin{prop} (\cite{ta} theorem 2) 
Let $A_1$ and $A_2$ be two unital $C(X)$-algebras. 

If $\|\;\|_\beta$ is a \mbox{${\rm C}^*$}-semi-norm on $A_1\otimes_{alg} A_2$ 
whose kernel is ${\cal J} (A_1,A_2)$, then 
\begin{center} 
$\forall\alpha\in A_1\otimes_{alg} A_2,\quad \|\alpha\|_\beta\geq
\|\alpha\|_m$.
\end{center}
\end{prop}{\mbox{{\bf Proof: }}} 
If we show that $S_\beta 
=P(A_1)\times_XP(A_2)$, then for every $\rho\in S_\beta$ and every $\alpha$ in 
$A_1\otimes_{alg} A_2$, we have $\rho (s^*\alpha^*\alpha s)\leq
\rho (s^*s)\|\alpha\|_\beta^2$ for all $s\in A_1\otimes_{alg} A_2$. Therefore  
\begin{center} 
\begin{tabbing}
\hspace{10pt}$\| \alpha\|_m{}^2$
\=$=\sup\Bigl\{\| (\sigma^x_1\otimes_{\min}\sigma_2^x)
(\alpha )\|^2 ,x\in X\Bigr\}$\\
\>$=\sup\Bigl\{ {\displaystyle
\frac{(\omega_1\otimes\omega_2)(s^*\alpha^*\alpha s)}
{(\omega_1\otimes\omega_2)(s^*s)}}$\=$, 
(\omega_1,\omega_2)\in P(A_1)\times_XP(A_2)$ and \\ 
\>\>$ s\in A_1\otimes_{alg} A_2\mbox{ such that } 
(\omega_1\otimes\omega_2)(s^*s)\not =0\Bigr\}$\\
\>$\leq\|\alpha\|_\beta{}^2$.
\end{tabbing}
\end{center}

\vskip12pt
Suppose that $S_\beta\not = P(A_1)\times_XP(A_2)$. 
Then there exist thanks to lemma \ref{elise} 
self-adjoint elements $a_i\in A_i$ and a point $x\in X$ such that 
$(a_1)_x\otimes (a_2)_x\not =0$ but ${(\omega_1\otimes\omega_2)
(a_1\otimes a_2)}=0$ for all couples $(\omega_1,\omega_2)\in S_\beta$. 

Let $B$ be the unital abelian $C(X)$-algebra generated by $C(X)$ and $a_1$ in 
$A_1$. The preceding lemma implies that 
$B{\mathop{\otimes}\limits^m}_{C(X)} A_2$ maps injectively into the 
Hausdorff completion $A_\beta$ of $A_1\otimes_{alg} A_2$ for $\|\;\|_\beta$. 

Consider pure states $\rho\in P(B_x)$ and $\omega_2\in P(\, (A_2)_x)$ such that 
$\rho (a_1)\not =0$ and $\omega_2(a_2)\not =0$ and extend the pure state 
$\rho\otimes\omega_2$ on $B{\mathop{\otimes}\limits^m}_{C(X)} A_2$ to a pure 
state $\omega$ on $A_\beta$. 
If we set $\omega_1(a)=\omega (a\otimes 1)$ for $a\in A_1$, then $\omega_1$ is 
pure and $\omega (\alpha )=(\omega_1\otimes\omega_2)(\alpha )$ for all 
$\alpha\in A_1\otimes_{alg} A_2$ since $\omega$ and $\omega_2$ are pure 
(\cite{ta} lemma 4). 
As a consequence, $(\omega_1,\omega_2)\in S_\beta$, which is absurd since 
$(\omega_1\otimes\omega_2)(a_1\otimes a_2)=
\rho (a_1)\omega_2(a_2)\not =0$. $\Box$

\begin{prop} Given two $C(X)$-algebras $A_1$ and $A_2$, 
the semi-norm $\|\,\|_m$ defines the minimal \mbox{${\rm C}^*$}-norm on the 
involutive algebra $(A_1\otimes_{alg} A_2)/{\cal J} (A_1,A_2)$.
\end{prop}{\mbox{{\bf Proof: }}} 
Let $\|\;\|_\beta$ be a \mbox{${\rm C}^*$}-norm on 
$(A_1\otimes_{alg} A_2)/{\cal J} (A_1,A_2)$. Thanks 
to the previous proposition, all we need to prove is that one can extend 
$\|\;\|_\beta$ to a \mbox{${\rm C}^*$}-norm on 
${({\cal A}_1\otimes_{alg}{\cal A}_2)}/{\cal J} ({\cal A}_1,{\cal A}_2)$, where 
${\cal A}_1$ and ${\cal A}_2$ are the unital $C(X)$-algebras associated 
to the $C(X)$-algebras $A_1$ and $A_2$ (definition \ref{lakral}). 

\vskip12pt
Consider the Hausdorff completion $D_\beta$ of 
$(A_1\otimes_{alg} A_2)/{\cal J} (A_1,A_2)$ 
and denote by $\pi_i$ the canonical representation of $A_i$ in $M(D_\beta )$ 
for $i=1,2$. Let us define the representation $\widetilde{\pi}_i$ of 
${\cal A}_i$ in $M(D_\beta\oplus A_1\oplus A_2\oplus C(X))$ by the 
following formulas: 
\begin{center}
$\begin{array}{rl}
\widetilde{\pi}_1(b_1+u(f))(\alpha\oplus a_1\oplus a_2\oplus g)&=
(\pi_1(b_1)+g)\alpha\oplus (b_1+f)a_1\oplus fa_2\oplus fg \\
\widetilde{\pi}_2(b_2+u(f))(\alpha\oplus a_1\oplus a_2\oplus g)&=
(\pi_2(b_2)+g)\alpha\oplus fa_1\oplus (b_2+f)a_2\oplus fg
\end{array}$
\end{center}
For $i=1,2$, let $\varepsilon_i:{\cal A}_i\rightarrow C(X)$ be the map defined 
by 
\begin{center} 
$\varepsilon_i[a+u(f)]=f$ for $a\in A_i$ and $f\in C(X)$. 
\end{center}
Then using the maps 
$(\varepsilon_1\otimes\varepsilon_2)$, $(\varepsilon_1\otimes id)$ and 
$(id\otimes\varepsilon_2)$, one proves 
easily that if $\alpha\in{\cal A}_1\otimes_{alg}{\cal A}_2$, 
$(\widetilde{\pi}_1\otimes\widetilde{\pi}_2)(\alpha )=0$ if and only if 
$\alpha$ belongs to ${\cal J} ({\cal A}_1,{\cal A}_2)$. 

Therefore, the norm of $M(D_\beta\oplus A_1\oplus A_2\oplus C(X))$ 
restricted to the subalgebra 
$({\cal A}_1\otimes_{alg}{\cal A}_2)/{\cal J} ({\cal A}_1,{\cal A}_2)$ 
extends $\|\;\|_\beta$. $\Box$

\vskip12pt\noindent {\sl Remark:} As the $C(X)$-algebra $A$ is nuclear if and 
only if every fibre $A_x$ is nuclear (\cite{la}), 
$A{\mathop{\otimes}\limits^M}_{C(X)} B\simeq 
A{\mathop{\otimes}\limits^m}_{C(X)} B$ for 
every $C(X)$-algebra $B$ if and only if $A$ is nuclear. 

\section{$\!\!\!$When does the equality ${\cal I} (A,B)\! =
\!\!{\cal J} (A,B)$ hold?}
\indent\indent Given two $C(X)$-algebras $A$ and $B$, 
Giordano and Mingo have studied in \cite{gior} the case where the algebra 
$C(X)$ is a von Neumann algebra: their theorem~3.1 and lemma~1.5 of \cite{kus} 
imply that in that case, we always have the equality ${\cal I} (A,B) =
{\cal J} (A,B)$. 

Our purpose in this section is to find sufficient conditions on the 
$C(X)$-algebras $A$ and $B$ in order to ensure this equality and to present a 
counter-example in the general case.

\begin{prop}\label{gogol} Let $X$ be a second countable Hausdorff compact space 
and let $A$ and $B$ be two $C(X)$-algebras.

If $A$ is a continuous field of \mbox{${\rm C}^*$}-algebras over $X$, then 
${\cal I} (A,B) ={\cal J} (A,B)$. 
\end{prop}{\mbox{{\bf Proof: }}} 
Let us prove by induction on the non negative integer $n$ that if 
\begin{center}
$s={\displaystyle \sum_{1\leq i\leq n}}a_i\otimes b_i\in{\cal J} (A,B)$, 
\end{center}
then $s$ belongs to the ideal ${\cal I} (A,B)$.

\vskip12pt If $n=0$, there is nothing to prove. Consider therefore an 
integer $n>0$ and suppose the result has been proved for any $p<n$.

Fix an element $s={\displaystyle {\sum_{1\leq i\leq n}}a_i\otimes 
b_i}\in{\cal J} (A,B)$ and define the continuous positive 
function $h\in C(X)$ by the formula $h(x)^{10}=\sum\| (a_k)_x\|^2$.

The element $a_k'=h^{-4}a_k$ is then well defined in $A$ for every $k$ since 
$a_k^*a_k\leq h^{10}$. Consequently, the function $f_k(x)=\| (a_k')_x\|$ is 
continuous. 

\vskip12pt 
For $1\leq k\leq n$, let $D_k$ denote the separable $C(X)$-algebra generated by 
$1$ and the $a_k'{}^*a_j'\,$,
$1\leq j\leq n$, in the unital $C(X)$-algebra ${\cal A}$ associated 
to $A$ (definition~\ref{lakral}). Then $D_k$ is a unital continuous field of 
\mbox{${\rm C}^*$}-algebras over $X$ (see for instance 
\cite{bla1} proposition 3.2).

\vskip12pt
Consider the open subset ${\cal S}^k=\{\psi\in {\cal S}_X(D_k)/\, 
\psi [a_k'{}^*a_k']>\psi (f_k^2/2)\}$. 
If we apply lemma 3.6~{\it b)} of \cite{bla1} to the 
restriction of $p:{\cal S}_{X} (D_k)\rightarrow X$ to ${\cal S}^k$, 
we may construct a continuous field of states $\omega_k$ on $D_k$ such 
that $\omega_k [a_k'{}^*a_k']\geq f_k^2/2$.

\vskip12pt 
Now, if we set $s'=\sum_i a_i'\otimes b_i$, as $(a_k'{}^*\otimes 1)s'$ belongs 
to ${\cal J} (D_k,B)$, 
\begin{center}
$(\omega_k\otimes id)[(a_k'{}^*\otimes 1)s']=\omega_k[a_k'{}^*a_k']b_k+
\mathop{\sum}\limits_{j\not =k}\omega_k[a_k'{}^*a_j']b_j=0$.
\end{center}
Noticing that $f^3_k$ is in the ideal of $C(X)$ generated by 
$\omega_k [a_k'{}^*a_k']$, we get that $f^3_kb_k$ belongs to the $C(X)$-module 
generated by the $b_j$, $j\not =k$, and thanks to the induction hypothesis, it 
follows that $(f_k^3\otimes 1)s'\in {\cal I} (A,B)$ for each $k$.

But $h^2=\sum_k f_k^2$ and so $h^4\leq n\sum_k f_k^4$ is in the ideal of $C(X)$ 
generated by the $f_k^3$, which implies 
$s=(h^4\otimes 1)s'\in{\cal I} (A,B)$. $\Box$

\vskip12pt\noindent {\sl Remarks:} 
{\bf 1.} As a matter of fact, it is not necessary to 
assume that the space $X$ satisfies the second axiom of countability thanks to 
the following lemma of \cite{kirch}: if $P(a)\in C(X)$ denotes the map 
$x\mapsto\| a_x\|$ for $a\in{\cal A}$, there exists a separable 
\mbox{${\rm C}^*$}-subalgebra 
$C(Y)$ of $C(X)$ with same unit such that if $D_k$ is the separable unital 
\mbox{${\rm C}^*$}-algebra generated by $C(Y).1\subset{\cal A}$ and the 
$a_k'{}^*a_j'\,$,
$1\leq j\leq n$, then $P(D_k)=C(Y)$. Furthermore, if $\Phi :X\rightarrow Y$ 
is the transpose of the inclusion map $C(Y)\hookrightarrow C(X)$ restricted to 
pure states, the map $D/(C_{\Phi (x)}(Y)D)\rightarrow A_x$ is a 
monomorphism for every $x\in X$ since ${\cal A}$ is continuous. 

Consider now a continuous field of states $\omega_k:D_k\rightarrow C(Y)$ on 
the continuous field $D_k$ over $Y$; 
if $\sum c^j\otimes d^j\in D_k\otimes_{alg} B$ is 
zero in $A_x\otimes_{alg} B_x$ for $x\in X$, then $\sum\omega_k(c^j)(x)d^j_x=0$ 
in $B_x$, which enables us to conclude as in the separable case.

\vskip12pt\noindent {\bf 2.} J. Mingo has drawn the author's attention to the 
following result of Glimm (\cite{glad} lemma 10): if $C(X)$ is a von Neumann 
algebra and $A$ is a $C(X)$-algebras, then 
$P(a)^2=\min\{ z\in C(X)_+, \, z\geq a^*a\}$ 
is continuous for every $a\in{\cal A}$. Therefore we always have in 
that case the equality ${\cal I} (A,B) ={\cal J} (A,B)$ thanks to the 
previous remark. 

\begin{cor}\label{bingo} Let $A$ and $B$ be two $C(X)$-algebras and 
assume that there exists a finite subset $F=\{ x_1,\cdots ,x_p\} 
\subset X$ such that for all $a\in A$, the function $x\mapsto\| a_x\|$ is 
continuous on $X\backslash F$.

Then the ideals ${\cal I} (A,B)$ and ${\cal J} (A,B)$ coincide.
\end{cor}{\mbox{{\bf Proof: }}}
Fix an element $\alpha =\sum_{1\leq i\leq n} a_i\otimes b_i\in 
A\otimes_{alg} B$ which belongs to ${\cal J} (A,B)$. In particular, we have 
$\sum_i (a_i)_{x}\otimes (b_i)_{x}=0$ in $A_x\otimes_{alg} B_x$ for each 
$x\in F$.

As a consequence, thanks to theorem III of \cite{mu}, we may find complex 
matrices $(\lambda_{i,j}^m)_{i,j}\in M_n(\C )$ for all $1\leq m\leq p$ such 
that, if we define the elements $c_k^m\in A$ and $d_k^m\in B$ by the formulas 
\begin{center} 
$c_k^m=\sum_{i}\lambda_{i,k}^ma_i$ and 
$d_k^m=b_k-\sum_j \lambda_{k,j}^mb_j$, 
\end{center}
we have $(c_k^m)_{x_m}=0$ and $(d_k^m)_{x_m}=0$ for all $k$ and all $m$.

Consider now a partition $\{ f_l\}_{1\leq l\leq p}$ of $1\in C(X)$ such that 
for all $1\leq l,m\leq p$, $f_l(x_m)=\delta_{l,m}$ where $\delta$ is the 
Kronecker symbol and define for all $1\leq k\leq n$ the elements 
$c_k=\sum_m f_m c_k^m$ and $d_k=\sum_m f_m d_k^m$. Thus, 
\begin{center}
$\begin{array}{rl}
\alpha&=\left(\sum_{i}a_i\otimes d_i\right) 
+\left(\sum_{i,j,m}\lambda_{i,j}^ma_i\otimes f_mb_j\right)\\
&=\left(\sum_{i}a_i\otimes d_i\right) +\left(\sum_{j}c_j\otimes b_j\right)+ 
\left(\sum_{i,j,m}\lambda_{i,j}^m (a_i\otimes f_mb_j-f_ma_i\otimes b_j)\right)
\end{array}$ 
\end{center} and 
there exists therefore an element $\beta\in{\cal I} (A,B)$ such that 
$\alpha -\beta$ 
admits a finite decomposition $\sum_i a_i'\otimes b_i'$ with 
$a_i'\in C_0(X\backslash F)A$ and $b_i'\in C_0(X\backslash F)B$.

But $C_0(X\backslash F)A$ is a continuous field. Accordingly, 
proposition \ref{gogol} implies that $\alpha -\beta\in
{\cal I} (C_0(X\backslash F)A,C_0(X\backslash F)B)\subset 
{\cal I} (A,B)$. $\Box$

\vskip12pt\noindent {\sl Remark:} If $\widetilde{\N}=\N\cup\{\infty\}$ is 
the Alexandroff compactification of $\N$ and if $A$ and $B$ are two 
$C(\widetilde{\N})$-algebras, the corollary \ref{bingo} implies the equality 
${\cal I} (A,B) ={\cal J} (A,B)$.

\vskip22pt
Let us now introduce a counter-example in the general case.

Consider a dense countable subset $X=\{ a_n\}_{n\in\N}$ of the interval 
$[0,1]$. 
The \mbox{${\rm C}^*$}-algebra $C_0(\N )$ of sequences with values in $\C$ 
vanishing at infinity is then endowed with the $C([0,1])$-algebra structure 
defined by: 
\begin{center}
$\forall f\in C([0,1])$, $\forall\alpha =(\alpha_n)\in C_0(\N )$, 
$\quad(f.\alpha )_n =f(a_n)\alpha_n$ for $n\in\N$.
\end{center}
If we call $A$ this $C([0,1])$-algebra, then $A_x=0$ for all $x\not\in X$. 

Indeed, assume that $x\not\in X$ and take $\alpha\in A$. If $\varepsilon >0$, 
there exists $N\in\N$ such that $|\alpha_n| <\varepsilon$ for all $n\geq N$. 
Consider a continuous function $f\in C([0,1])$ such that $0\leq f\leq 1$, 
$f(x)=0$ and $f(a_i)=1$ for every $1\leq i\leq N$; we then have: 
\begin{center}
$\|\alpha_x\|\leq\| (1-f)\alpha\| <\varepsilon$.
\end{center}

Let $Y=\{ b_n\}_{n\in\N}$ be another dense countable subset of $[0,1]$ and 
denote by $B$ the associated $C([0,1])$-algebra whose underlying algebra 
is $C_0(\N )$. 

Then, if $X\cap Y=\emptyset$, the previous remark implies that for 
every $x\in [0,1]$, $A_x\otimes_{alg} B_x=0$ and hence, 
${\cal J} (A,B) =A\otimes_{alg} B$. 
What we therefore need to prove is that the ideal ${\cal I} (A,B)$ is strictly 
included in $A\otimes_{alg} B$.

\vskip12pt
Let us fix two sequences $\alpha\in A$ and $\beta\in B$ whose terms are all 
non zero and suppose that $\alpha\otimes\beta$ admits a decomposition 
$\sum_{1\leq i\leq k}[(f_i\alpha^i)\otimes\beta^i -
\alpha^i\otimes (f_i\beta^i)]$ in $A\otimes_{alg} B$. Then for every 
$n,m\in\N$,
$$\alpha_n\beta_m =\sum_{1\leq i\leq k}
\alpha^i_n\beta^i_m[f_i(a_n)-f_i(b_m)].$$ 
Now, if we set $\phi_i(a_n) = \alpha_n^i/\alpha_n$ and $\psi_i(b_n) 
=\beta_n^i/\beta_n$ for $1\leq i\leq k$ and $n\in\N$, this equality means that 
for all $(x,y)\in X\times Y$, 
\begin{center}
$1={\sum}_{1\leq i\leq k}\psi_i (x)\varphi_i(y)[f_i(x)-f_i(y)]$. 
\end{center}
But this is impossible because of the following proposition:

\begin{prop}
Let $X$ and $Y$ be two dense subsets of the interval $[0,1]$ and let $n$ be a 
non negative integer. 

Given continuous functions $f_i$ on $[0,1]$ and numerical functions 
$\psi_i:X\rightarrow\C$ and $\varphi_i:Y\rightarrow\C$ for $1\leq i\leq n$, if 
there exists a constant $c\in\C$ such that 
\begin{center}
 $\forall (x,y)\in X\times Y$, 
$\quad \mathop{\sum}\limits_{1\leq i\leq n}\psi_i 
(x)\varphi_i(y)[f_i(x)-f_i(y)]=c$ 
\end{center}
then $c=0$.
\end{prop}{\mbox{{\bf Proof: }}} 
We shall prove the proposition by induction on $n$.

\vskip12pt
If $n=0$, the result is trivial. Take therefore $n>0$ and assume that the 
proposition is true for any $k<n$. 

Suppose then that the subsets $X$ and $Y$ of $[0,1]$, the 
functions $f_i$, $\psi_i$ and $\varphi_i$, $1\leq i\leq n$, satisfy the 
hypothesis of the proposition for the constant $c$. 

\vskip12pt
For $x\in X$, let $p(x)\leq n$ be the dimension of the vector space 
generated in $\C^n$ by the $\Bigl(\varphi_i(y)[f_i(x)-f_i(y)]
\Bigr)_{1\leq i\leq n}$, $y\in Y$. 

If $p(x)<n$, there exists a subset $F(x)\subset\{ 1,\ldots ,n\}$ of cardinal 
$p(x)$ such that for every $j\not\in F(x)$:
\begin{center}
$\varphi_j(y)[f_j(x)-f_j(y)]=\sum_{i\in F(x)} \lambda^j_i(x) \varphi_i(y) 
[f_i(x)-f_i(y)]$ for all $y\in Y$, 
\end{center}
where the $\lambda^j_i(x)\in\C$ are given by the Cramer formulas. As a 
consequence, 
\begin{center} $\sum_{i\in F(x)}\biggl( \psi_i(x)+\sum_{j\not\in F(x)}
[\lambda_i^j(x)\psi_j(x)]\biggr)\varphi_i(y)[f_i(x)-f_j(y)]=c$.
\end{center}

Now, if $p(x)<n$ for every $x\in X$, there exists a subset $F\subset 
\{ 1,\ldots ,n\}$ of cardinal $p<n$ such that the interior of the closure of 
the set of those $x$ for which $F(x)=F$ is not empty and contains therefore 
a closed interval homeomorphic to $[0,1]$. The induction hypothesis for $k=p$ 
implies that $c=0$.

\vskip12pt
Assume on the other hand that $x_0\in X$ verifies $p(x_0)=n$. 
We may then find $y_1, \cdots y_n$ in $Y$ such that if we set 
\begin{center} 
$a_{i,j}(x)= \varphi_i(y_j)[f_i(x)-f_i(y_j)]$, 
\end{center}
the matrix $\left( a_{i,j}(x_0)\right)$ is invertible. There exists therefore a 
closed connected neighborhood $I$ of $x_0$ on which the matrix 
$\left( a_{i,j}(x)\right)$ remains invertible.

But for each $1\leq j\leq n$, $\mathop{\sum}\limits_{i} a_{i,j}(x)\psi_i(x)=c$ 
and therefore the $\psi_i(x)$ extend by the Cramer 
formulas to continuous functions on the closed interval $I$.

\vskip12pt
For $y\in Y\cap I$, let $q(y)$ denote the dimension of the vector space 
generated in $\C^n$ by the $\Bigl(\psi_i (x)[f_i(x)-f_i(y)]
\Bigr)_{1\leq i\leq n}$, $x\in X\cap I$. 

If $q(y)<n$ for every $y$, then the induction hypothesis implies $c=0$. But if 
there exists $y_0$ such that $q(y_0)=n$, we may find an interval $J\subset I$ 
homeomorphic to $[0,1]$ on which the $\varphi_i$ extend to continuous 
functions; evaluating the starting formula at a point $(x,x)\in J\times J$, 
we get $c=0$. $\Box$

\section{The associativity}
\indent\indent 
Given three $C(X)$-algebras $A_1$, $A_2$ and $A_3$, we deduce from 
\cite{bla1} corollaire 3.17: 
\begin{center}
$[(A_1{\mathop{\otimes}\limits^M}_{C(X)} A_2){\mathop{\otimes}\limits^M}_{C(X)} 
A_3]_x\! =\! (A_1{\mathop{\otimes}\limits^M}_{C(X)} A_2)_x\otimes_{\max}\! 
(A_3)_x\! =\! (A_1)_x\otimes_{\max}\! (A_2)_x\otimes_{\max}\! (A_3)_x$, 
\end{center}
which implies the associativity of the tensor product 
$\cdot{\mathop{\otimes}\limits^M}_{C(X)}\cdot$ over $C(X)$.

\vskip12pt
On the contrary, the minimal tensor product 
$\cdot{\mathop{\otimes}\limits^m}_{C(X)}\cdot$ over $C(X)$ is not 
in general associative. 
Indeed, Kirchberg and Wassermann have shown in \cite{kirch} that if 
$\widetilde{\N} =\N\cup\{\infty\}$ is the Alexandroff compactification of 
$\N$, there exist separable continuous fields $A$ and $B$ such that 
\begin{center}
$(A{\mathop{\otimes}\limits^m}_{C(\widetilde{\N})}B)_\infty\not =
A_\infty\otimes_{\min} B_\infty$.
\end{center}
If we now endow the \mbox{${\rm C}^*$}-algebra $D=\C$ with the 
$C(\widetilde{\N })$-algebra structure defined by $f.a=f(\infty )a$, then for 
all $C(\widetilde{\N })$-algebra $D'$, we have 
\begin{center}
$[D{\mathop{\otimes}\limits^m}_{C(\widetilde{\N})} D']_n=
\left\{\begin{array}{cl} 
0&\mbox{ if $n$ is finite,}\\
(D')_\infty&\mbox{ if $n=\infty$.}
\end{array}\right.$
\end{center} 
Therefore, ${[(A{\mathop{\otimes}\limits^m}_{C(\widetilde{\N})} B)
{\mathop{\otimes}\limits^m}_{C(\widetilde{\N})}D]_\infty}\simeq 
{(A{\mathop{\otimes}\limits^m}_{C(\widetilde{\N})} B)_\infty}$ whereas
${[A{\mathop{\otimes}\limits^m}_{C(\widetilde{\N})} 
(B{\mathop{\otimes}\limits^m}_{C(\widetilde{\N})} D)]_\infty}$ is 
isomorphic to ${A_\infty\otimes_{\min} B_\infty}$.

\vskip12pt
However, in the case of (separable) continuous fields, we can deduce the 
associativity of $\cdot{\mathop{\otimes}\limits^m}_{C(X)}\cdot$ from the 
following proposition: 
\begin{prop} Let $A$ and $B$ be two $C(X)$-algebras.

Assume $\pi$ is a field of faithful representations of $A$ in the 
Hilbert $C(X)$-module ${\cal E}$, then the morphism 
$a\otimes b\rightarrow \pi (a)\otimes b$ 
induces a faithful $C(X)$-linear representation of 
$A{\mathop{\otimes}\limits^m}_{C(X)} B$ in the Hilbert 
$B$-module ${\cal E}\otimes_{C(X)} B$. 
\end{prop}{\mbox{{\bf Proof: }}} 
Notice that for all $x\in X$, we have 
$({\cal E}\otimes_{C(X)} B)\otimes_B B_x={\cal E}_x\otimes B_x$. 

Now, as $B$ maps injectively in $B_d=\oplus_{x\in X}B_x$, 
${\cal L}_B({\cal E}\otimes_{C(X)} B)$ maps injectively in $\oplus_{x\in X} 
{\cal L}_{B_x}({\cal E}_x\otimes B_x)\subset
{\cal L}_{B_d}({\cal E}\otimes_{C(X)} B\otimes_{B}B_d)$ and 
therefore if $\alpha\in A\otimes_{alg} B$, we have 
$\| (\pi\otimes id)(\alpha )\| =\mathop{\sup}\limits_{x\in X}\| 
(\pi_x\otimes id)(\alpha_x)\| =\|\alpha\|_m$. $\Box$

\vskip12pt
Accordingly, if for $1\leq i\leq 3$, $A_i$ is a separable continuous field of 
\mbox{${\rm C}^*$}-algebras over $X$ which admits a field of faithful 
representations in the $C(X)$-module ${\cal E}_i$, the $C(X)$-representations 
of ${(A_1{\mathop{\otimes}\limits^m}_{C(X)} A_2)
{\mathop{\otimes}\limits^m}_{C(X)} A_3}$ and 
${A_1{\mathop{\otimes}\limits^m}_{C(X)} 
(A_2{\mathop{\otimes}\limits^m}_{C(X)} A_3)}$ in the Hilbert $C(X)$-module 
${({\cal E}_1\otimes_{C(X)}{\cal E}_2)\otimes_{C(X)}{\cal E}_3}=
{{\cal E}_1\otimes_{C(X)} ({\cal E}_2\otimes_{C(X)}{\cal E}_3)}$ are faithful, 
and hence the maps 
$A_1\otimes_{\min} A_2\otimes_{\min} A_3\rightarrow (A_1{\mathop{\otimes}\limits^m}_{C(X)} A_2)
{\mathop{\otimes}\limits^m}_{C(X)} A_3$ and 
$A_1\otimes_{\min} A_2\otimes_{\min} A_3\rightarrow 
A_1{\mathop{\otimes}\limits^m}_{C(X)} 
(A_2{\mathop{\otimes}\limits^m}_{C(X)} A_3)$ have the same kernel.

\vskip12pt
{\small {\it 
Etienne Blanchard,\ \ Mathematisches Institut, Universit\"at Heidelberg\\ 
Im Neuenheimer Feld 288,\ \ D $-$ 69120 Heidelberg\\
e-mail: blanchar@mathi.uni-heidelberg.de }} 

\end{document}